\documentclass[12pt]{article}
\hoffset - 10 mm

\catcode`\@=11 \@addtoreset{equation}{section}

\catcode`\@=12

\newtheorem{Theorem}{Theorem}[section]
\newtheorem{Proposition}{Proposition}[section]
\newtheorem{Lemma}{Lemma}[section]
\newtheorem{Corollary}{Corollary}[section]

\newcommand{\bTheorem}[1]{
\begin{Theorem} \label{T#1} }
\newcommand{\eT}{\end{Theorem}}

\newcommand{\bProposition}[1]{
\begin{Proposition} \label{P#1}}
\newcommand{\eP}{\end{Proposition}}

\newcommand{\bLemma}[1]{
\begin{Lemma} \label{L#1} }
\newcommand{\eL}{\end{Lemma}}

\newcommand{\bCorollary}[1]{
\begin{Corollary} \label{C#1} }
\newcommand{\eC}{\end{Corollary}}

\newcommand{\bFormula}[1]{
\begin{equation} \label{#1}}
\newcommand{\eF}{\end{equation}}
\newcommand{\Ome}{\Omega_\ep}
\newcommand{\Den}{\Delta_{\ep,N}}
\newcommand{\Ov}[1]{\overline{#1}}
\newcommand{\DC}{C^\infty_c}

\newcommand{\vr}{\varrho}
\newcommand{\vre}{\vr_\ep}

\newcommand{\vue}{\vu_\ep}

\newcommand{\vu}{\vc{u}}
\newcommand{\vc}[1]{{\bf #1}}
\newcommand{\qed}{\bigskip \rightline {Q.E.D.} \bigskip}
\newcommand{\Div}{{\rm div}_x}
\newcommand{\Grad}{\nabla_x}

\newcommand{\tn}[1]{\mbox {\F #1}}
\newcommand{\dx}{{\rm d} {x}}
\newcommand{\dt}{{\rm d} t }

\newcommand{\dxdt}{\dx \ \dt}
\newcommand{\intO}[1]{\int_{\Omega} #1 \ \dx}
\newcommand{\intOe}[1]{\int_{\Omega_\ep} #1 \ \dx}

\newcommand{\bProof}{{\bf Proof: }}

\newcommand{\ep}{\varepsilon}

\font\F=msbm10 scaled 1000
\newcommand{\R}{\mbox{\F R}}

\date{}
\begin{document}

\title{Stability with respect to domain of the low Mach number limit of
compressible viscous fluids}
\author{Eduard Feireisl\thanks{The work of E.F.
was supported by Grant 201/09/0917 of GA \v CR as a part of the
general research programme of the Academy of Sciences of the Czech
Republic, Institutional Research Plan AV0Z10190503.} \and Trygve Karper
\thanks{The work of T.K. was supported by Ne\v cas Center of Mathematical Modelling (LC06052).} \and
Ond\v rej Kreml\thanks{The work of O.K. is a part of  the general
research programme of the Academy of Sciences of the Czech Republic,
Institutional Research Plan AV0Z10190503.} \and Jan Stebel
\thanks{The work of J.S. was supported by Grant 201/09/0917 of GA \v
CR.}} \maketitle

\centerline{Institute of Mathematics, Academy of Sciences of the Czech Republic}
\centerline{\v Zitn\' a 25, 115 67 Praha 1, Czech Republic}

\centerline{and}
\centerline{$^\dagger$CSCAMM, University of Maryland}
\centerline{4146 CSIC Building \#406,
Paint Branch Drive,
College Park MD 20742, USA}

\begin{abstract}

We study the asymptotic limit of solutions to the barotropic
Navier-Stokes system, when the Mach number is proportional to a
small parameter $\ep \to 0$ and the fluid is confined to an exterior
spatial domain $\Omega_\ep$ that may vary with $\ep$. As $\epsilon \rightarrow 0$, it is shown that the
fluid density becomes constant while the velocity converges
 to a solenoidal vector field satisfying the incompressible
Navier-Stokes equations on a limit domain. The velocities approach
the limit strongly (a.a.) on any compact set, uniformly with respect
to a certain class of domains. The proof is based on spectral
analysis of the associated wave propagator (Neumann Laplacian)
governing the motion of acoustic waves.

\end{abstract}


\section{Introduction}
\label{i}

There is a vast number of mathematical models of
\emph{incompressible} fluids that can be identified as a singular
limit of more complex systems describing the motion of compressible
and/or heat conducting fluids, see the review papers by Gallagher
\cite{Gallag}, Masmoudi \cite{MAS}, \cite{MAS1}, or the monograph
\cite{FEINOV} and the references therein. In many cases the
resulting system is driven by an exterior force related to
gravitation of rigid objects \emph{outside} the fluid domain
whereas, at the same time, the fluid occupies the whole physical
space $R^3$, see the study of the Oberbeck-Boussinesq approximation
by Brandolese and Schonbek \cite{BranSchon}. Such a situation, if
physically relevant, can be viewed as a singular limit, where the
primitive system is posed on a family of \emph{exterior} domains
$\Omega_\ep$, with $R^3 \setminus \Omega_\ep$ being the rigid
body(ies) acting on the fluid by their gravitation, and $\Omega_\ep
\to R^3$ as $\ep \to 0$ in a certain sense. Our goal in the present
study is to develop a method for studying the incompressible limits,
with the Mach number ${\rm Ma} = \ep \to 0$, on a family of
\emph{exterior} domains $\Omega_\ep \subset R^3$ varying with $\ep >
0$. In particular, we identify a class of domains giving rise to
\emph{uniform} convergence of the fluid velocities, independent of
the specific shape of their boundaries. To this end, we adapt the
technique introduced in \cite{EF84}, based on spectral theory for
the corresponding acoustic wave propagator - the Neumann Laplacian
on $\Ome$. The adaptation leans on delicate estimates of the
associated Hemholtz projections based on the results by Farwig,
Kozono, and Sohr in \cite{FAKOSO} and \cite{FAKOSO1}.

For the sake of simplicity, we focus only on the mechanical aspects
of the fluid motion, ignoring completely the effect of temperature
changes. Accordingly, we consider the compressible
\emph{Navier-Stokes system} 
in
Eulerian
reference coordinates:
\bFormula{i1}
\partial_t \vr + \Div (\vr \vu) = 0,
\eF
\bFormula{i2}
\partial_t (\vr \vu) + \Div (\vr \vu \otimes \vu) + \frac{1}{\ep^2} \Grad p(\vr) = \Div \tn{S}(\Grad \vu),
\eF
\bFormula{i3}
\tn{S}(\Grad \vu) = \mu \left( \Grad \vu + \Grad^t \vu - \frac{2}{3} \Div \vu \tn{I} \right), \ \mu > 0,
\eF
where the unknown functions are the mass density $\vr = \vr(t,x)$ and
the vector velocity $\vu = \vu(t,x)$, with $t \in (0,T)$, and $x \in \Omega_\ep$. The symbol $p = p(\vr)$ denotes the pressure - a given function of the density - whereas $\tn{S}$ stands for the viscous stress tensor.
The small parameter $\ep$ represents the Mach number tending to zero
in the asymptotic limit.

The system (\ref{i1} - \ref{i3}) is supplemented by the (acoustically
hard) \emph{complete slip} boundary conditions
\bFormula{i4}
\vu \cdot \vc{n}|_{\partial \Omega_\ep} = 0, \ [\tn{S} \cdot \vc{n}] \times \vc{n}|_{\partial \Omega_\ep} = 0,
\eF
where the  $\vc{n}$ denotes the (outer) normal vector to
$\partial \Omega_\ep$. Moreover, as the fluid occupies an exterior
domain, the behavior of $\varrho$, $\vu$ at infinity must be
specified:
\bFormula{i4a}
\vr \to \Ov{\vr} > 0,\  \vu \to 0 \ \mbox{as}\ |x| \to \infty,
\eF
where $\Ov{\vr}$ is a constant.

In the singular limit $\ep \to 0$, the fluid is driven to
\emph{incompressibility} as the speed of sound tends to infinity.
The incompressible limit for system (\ref{i1} - \ref{i3}) as well as
related problems have been studied by many authors, see Alazard
\cite{AL2}, Klainerman and Majda \cite{KM2}, Lions and Masmoudi
\cite{LIMA1}, or Schochet \cite{SCH2}, to name only a few. Similarly
to Lions and Masmoudi \cite{LIMA1}, \cite{LIMA6}, our approach is
based on the concept of \emph{weak solutions} to the Navier-Stokes
system, where the presence of viscosity plays a crucial role. As is
well known (cf. \cite{LIMA1}), convergence of solutions of the system
(\ref{i1} - \ref{i3}) to the incompressible limit may be disturbed
by the presence of rapidly oscillating acoustic waves, here
represented by the gradient part of the velocity field. Since the
physical domains $\Omega_\ep$ are \emph{unbounded}, however, we
expect that acoustic waves disperse leaving very fast any bounded
part of the physical space as was observed by Alazard \cite{AL2}, Bresch
and Metivier \cite{BreMet}, Isozaki \cite{Isoz}, and \cite{FNP4}.
The main novelty of this paper is the fact that the physical domains
$\Omega_\ep$ are allowed to change their shape together with the
Mach number. In particular, we show that the rate of convergence is
uniform within a certain class of domains specified in Section
\ref{dc} below. An interesting aspect of the problem is the boundary
behavior of the limit velocity field $\vc{U}$. As shown in
\cite{BFNW}, the slip boundary conditions (\ref{i4}) may give rise
to the more standard \emph{no-slip} condition
\bFormula{i5}
\vc{U}|_{\partial \Omega} = 0,
\eF
or to the kind of
friction-driven boundary conditions identified in
\cite{BuFeNe}.

In order to see the principal difficulties involved, we rewrite the Navier-Stokes system in the form of \emph{Lighthill's acoustic analogy} \cite{LIGHTHILL1}, \cite{LIGHTHILL2}:
\bFormula{i6}
\ep \partial_t r + \Div \vc{V} = 0 \ \mbox{in} \ (0,T) \times \Omega_\ep,
\eF
\bFormula{i7}
\ep \partial_t \vc{V} + p'(\Ov{\vr}) \Grad r = \ep \Div \tn{L} \ \mbox{in}\
(0,T) \times \Omega_\ep,
\eF
supplemented with the boundary condition
\bFormula{i8}
\vc{V} \cdot \vc{n}|_{\partial \Omega_\ep} = 0,
\eF
where
\bFormula{i9}
r \equiv \frac{\vr - \Ov{\vr}}{\ep} , \ \vc{V} \equiv \vr \vu,
\eF
and $\tn{L}$ is the so-called Lighthill's tensor,
\bFormula{i10}
\tn{L} \equiv \tn{S} - \vr \vu \otimes \vu - \frac{1}{\ep^2}
\Big( p (\vr) - p'(\Ov{\vr}) (\vr - \Ov{\vr}) - p(\Ov{\vr}) \Big) \tn{I}.
\eF
Applying, formally, the Helmholtz projection to (\ref{i7}) we obtain
the wave (acoustic) equation
\bFormula{i11}
\ep \partial_t r + \Delta \Phi = 0,
\eF
\bFormula{i12}
\ep \partial_t \Phi + p'(\Ov{\vr}) r = \ep \Delta^{-1}_{\ep,N} \Div
\Div \tn{L},
\eF
\bFormula{i13}
\Grad \Phi \cdot \vc{n}|_{\partial \Omega_\ep} = 0,
\eF
for the \emph{acoustic potential} $\Phi = \Delta^{-1}_{\ep,N} \Div
\vc{V}$, where the symbol $\Delta_{\ep,N}$ denotes the standard
Laplace operator supplemented with  homogeneous Neumann boundary
conditions on $\partial \Omega_\ep$.

Accordingly, our main task consists in:
\begin{itemize}
\item estimating the forcing term $\Delta^{-1}_{\ep,N} \Div \Div \tn{L}$
as well as the initial data in terms of a suitable power of the
operator $\Delta_{\ep,N}$, cf. \cite{FNP4};
\item evaluating local rate of decay of solutions to the wave equation (\ref{i11}), (\ref{i12}).
\end{itemize}
In general,  bounds of $\Delta^{-1}_{\ep,N} \Div \Div \tn{L}$ in
terms of $\Delta_{\ep,N}$ depend on the shape of $\partial
\Omega_\ep$, where the latter must be at least of class $C^{1,1}$ to
recover the standard $W^{2,p}-$theory, not available on less smooth,
say, Lipschitz domains, see Grisvard \cite{Gris}. On the other hand,
our method is applicable to families $\{ \Ome \}_{\ep >0}$,
whose smoothness parameters blow up for $\ep \to 0$. In particular,
they may approach a less regular domain in the asymptotic limit
and/or their boundaries may oscillate similarly to \cite{BFNW}, see
Section \ref{dc}. As a result, the forcing term involving
Lighthill's tensor become unbounded for $\ep \to 0$, and this defect
must be compensated by uniform dispersive estimates of order
$\sqrt{\ep}$ which we achieve in the spirit of a result
due to Kato \cite{Kato}.

It is interesting to note that similar results 
on
\emph{bounded} domains, supplemented with the no-slip boundary
condition (\ref{i5}), where convergence of the velocities is
enforced by a viscous boundary layer (see Desjardins et al.
\cite{DGLM}) seem very sensitive and much less stable with respect
to domain perturbations, in particular they completely fail  on
balls.

The organization of the paper is as follows. In Section \ref{m}, we
recall some known facts concerning the compressible Navier-Stokes
system (\ref{i1} - \ref{i3}), including the available existence
theory, introduce the principal hypotheses concerning the admissible
class of domains, and state our main result. Section \ref{ub}
reviews the standard uniform bounds on solutions to (\ref{i1} -
\ref{i3}) in $(0,T) \times \Omega_\ep$ independent of the singular
parameter $\ep \to 0$. In Section \ref{ae}, we introduce the
acoustic equation and identify the terms in Lighthill's tensor. In
particular, we deduce estimates on Lighthill's tensor in terms of
$\ep$, based on the $W^{2,p}-$theory for the Neumann Laplacian
$\Delta_{\ep, N}$ defined on $\Ome$. Section \ref{nl} deals with the
spectral theory of the operator $-\Delta_{\ep,N}$ in
$L^2(\Omega_\ep)$. We introduce the associated spectral measures and
show local decay of acoustic waves of rate $\sqrt{\ep}$. Finally,
the incompressible limit is performed in Section \ref{co}. In
addition, we shortly discuss the limit boundary conditions in the
spirit of \cite{BuFeNe}.

\section{Preliminaries and main result}
\label{m}

Throughout the paper, the pressure $p$ is a continuously
differentiable function of the density such that
\bFormula{m1}
p \in C[0,\infty) \cap C^1(0, \infty), \ p'(\vr) > 0 \ \mbox{for all}\ \vr > 0,
\lim_{\vr \to \infty} \frac{p'(\vr)}{\vr^{\gamma - 1}} = p_\infty > 0,
\eF
for a certain $\gamma > 3/2$.

Multiplying the momentum equation (\ref{i2}) by $\vu$ and integrating by
parts leads to the energy inequality:
\bFormula{m3}
\frac{{\rm d}}{{\rm d}t} \intOe{ E_\ep (\vr, \vu) (\tau, \cdot) } +
\intOe{ \tn{S}(\Grad \vu) : \Grad \vu } \leq 0, \ \tau \in (0,T),
\eF
where we have set
\[
E_\ep(\vr, \vu) \equiv  \frac{1}{2} \vr |\vu|^2 + \frac{1}{\ep^2}
\left( P(\vr) - P'(\Ov{\vr})(\vr - \Ov{\vr}) - P(\Ov{\vr}) \right),
\]
\bFormula{m4}
P(\vr) \equiv \vr \int_1^\vr \frac{p(z)}{z^2} \ {\rm d}z.
\eF

As $P''(\vr) = p'(\vr)/\vr > 0$, the function $P$ is strictly convex and
$P(\vr) - P'(\Ov{\vr})(\vr - \Ov{\vr}) - P(\Ov{\vr}) \approx c (\vr - \Ov{\vr})^2$
provided $\vr \approx \Ov{\vr}$. Consequently, with \emph{initial data} of the form
\bFormula{m5}
\vr(0, \cdot) = \vr_{0, \ep} = \Ov{\vr} + \ep r_{0, \ep},\
\vu(0, \cdot) = \vu_{0,\ep},
\eF
with
\bFormula{m6}
\| r_{0,\ep} \|_{L^2(\Ome)} + \| r_{0,\ep} \|_{L^\infty(\Ome)} +
\| \vu_{0, \ep} \|_{L^2(\Ome;R^3)} \leq c,
\eF
we get the total (mechanical) energy associated to the initial data
\[
\intOe{ E_\ep (\vr_{0,\ep}, \vu_{0,\ep} ) } = \intOe{ \left(
\frac{1}{2} \vr_{0,\ep} |\vu_{0,\ep}|^2 + \frac{1}{\ep^2} \left(
P(\vr_{0,\ep}) - P'(\Ov{\vr})(\vr_{0,\ep} - \Ov{\vr}) - P(\Ov{\vr})
\right) \right)},
\]
bounded uniformly for $\ep \to 0$.

\subsection{Weak solutions}

\label{ws}

We say that a pair of functions $\vr$, $\vu$ represents a \emph{weak
solution} to the Navier-Stokes system (\ref{i1} - \ref{i3}), with
the boundary conditions (\ref{i4}), (\ref{i4a}), and the initial
data (\ref{m5}) if:

\begin{itemize}
\item $\vr \geq 0$, $( \vr - \Ov{\vr} ) \in (L^2 + L^\gamma)(\Ome)$, $\vu
\in L^2(0,T; W^{1,2}(\Ome;R^3))$ such that $\vu \cdot \vc{n}|_{\partial \Ome} = 0$;
\item the equation of continuity (\ref{i1}) is satisfied in the sense of renormalized solutions:
\bFormula{m7}
\int_0^T \int_{\Ome} \Big[ \left( \vr + b(\vr) \right) \partial_t \varphi +
\left( \vr + b(\vr) \right) \vu \cdot \Grad \varphi
\eF
\[
 + \left(
b(\vr) - b'(\vr) \vr \right) \Div \vu \varphi \Big] \ \dxdt = -
\intOe{ \left( \vr_{0,\ep} + b(\vr_{0, \ep}) \right) \varphi (0,
\cdot) },
\]
for any test function $\varphi \in \DC([0,T) \times \Ov{\Omega}_\ep)$, and any
$b \in C^\infty[0,\infty)$, $b' \in \DC[0, \infty)$;
\item the momentum equation (\ref{i2}) holds in the sense of the integral identity
\bFormula{m8}
\int_0^T \intOe{ \Big[ \vr \vu \cdot \partial_t \varphi + \vr \vu \otimes \vu :
\Grad \varphi + \frac{1}{\ep^2} p(\vr) \Div \varphi \Big] } \ \dt
\eF
\[
=
\int_0^T \intOe{ \tn{S}(\Grad \vu) : \Grad \varphi } \ \dt - \intOe{
\vr_{0,\ep} \vu_{0,\ep} \cdot \varphi (0, \cdot) },
\]
 for any test function $\varphi \in \DC([0,T) \times \Ov{\Omega}_\ep; R^3)$,
$\Grad \varphi \cdot \vc{n}|_{\partial \Ome} = 0$;
\item the energy inequality
\bFormula{m9}
- \int_0^T E_\ep(\vr, \vu) \partial_t \psi \ \dt + \int_0^T \psi \intOe{ \tn{S}(\Grad \vu) :
\Grad \vu } \ \dt
\eF
\[
\leq \intOe{ E_\ep( \vr_{0,\ep}, \vu_{0,\ep}) \psi(0) },
\]
holds for any $\psi \in \DC[0,T)$, $\psi \geq 0$.

\end{itemize}

\subsection{Admissible domains}

\label{dc}

Motivated by \cite{BFNW}, we introduce a class of admissible domains
that allows for ``oscillating'' boundaries. As first observed by
Casado-Diaz, Fernandez-Cara, and Simon \cite{CFS}, such a family of
domains may give rise to the \emph{no-slip} boundary condition
(\ref{i5}) for the limit velocity field. A general description of
all possible limit boundary conditions was obtained in
\cite{BuFeNe}.

Introducing a cone
\[
C(x, \omega, \delta,\xi)=\{y \in R^3 \ | \ 0< |y-x | \le \delta,\
(y-x) \cdot \xi
> \cos (\omega) |y-x | \},
\]
with vertex at $x$, aperture $2\omega < \pi$, height $\delta$, and
orientation given by a unit vector $\xi$, we say that $\Omega_\ep$
satisfies the {uniform} $\delta-$cone condition if for any $x_0 \in
\partial \Ome$, there exists a unit vector $\xi _{x_0} \in \R^3$
such that
\[
C(x, \omega, \delta,\xi_{x_0} ) \subset {\Omega}_\ep \
\mbox{whenever}\ x \in {\Omega}_\ep ,\ |x - x_0| < \delta,
\]
see Henrot and Pierre \cite[Definition 2.4.1]{hepi05}.

Assume we are given a family of domains $\{ \Ome \}_{\ep > 0}$
complying with the following hypotheses:

\begin{itemize}
\item
$\Ome \subset R^3$ is an exterior domain with $C^2-$boundary for
each fixed $\ep > 0$;
\item
there is a $d > 0$ such that
\[
R^3 \setminus \Ome \subset B_d \equiv \{ x \in R^3 \ | \ |x| < d \}
\ \mbox{for all}\ \ep > 0;
\]
\item
$\Omega_\ep$ satisfy the uniform $\delta-$cone condition with  
$\delta > 0$ (and $\omega$) independent of $\ep$;
\item for each $x_0 \in \partial \Ome$, there are two (open) balls
$B_r[x_i] \equiv \{ x ; |x-x_i | < r \}\subset \Ome$, $B_r[x_e]
\subset R^3 \setminus \Ome$ of radius $r > c_b \ep^\beta$  such that
\bFormula{dc1}
\Ov{B_r[x_i]} \cap \Ov{B_r[x_e]} = x_0,
\eF
with $c_b > 0$, $\beta > 0$ independent of $\ep$.

\end{itemize}

The above hypotheses give rise to the following properties enjoyed
by the family $\{ \Ome \}_{\ep > 0}$:

\begin{itemize}

\item \emph{Uniform extension property (see Jones \cite{Jones}).}

There exists an extension operator $E_\ep$,
\bFormula{uep}
E_\ep : W^{1,p}(\Ome) \to W^{1,p}(R^3),\
E_\ep [v] |_{\Ome} = v, \ \| E_\ep [v] \|_{W^{1,p}(R^3)} \leq c \| v \|_{W^{1,p}(\Ome)},
\eF
where the constant $c$ is independent of $\ep \to 0$.

\item \emph{Uniform Korn's inequality (see \cite[Proposition 4.1] {BuFe1}).}

Let $\vc{v} \in W^{1,2}(\Ome \cap B; R^3)$, and $M \subset \Ome \cap
B$ such that $|M| > m > 0$, where $B$ is a bounded ball. Then
\bFormula{uki}
\| \vc{v} \|_{W^{1,2}(\Ome \cap B;R^3)}^2 \leq c(m) \left( \left\|
\Grad \vc{v} + \Grad^t \vc{v} - \frac{2}{3} \Div \vc{v} \tn{I}
\right\|_{L^2(\Ome \cap B; R^{3 \times 3})}^2 + \int_{M} |\vc{v} |^2
\ \dx \right),
\eF
with $c(m)$ independent of $\ep \to 0$.

\item \emph{Compactness (see Henrot and Pierre \cite[Theorem
2.4.10]{hepi05}).}

There exists an exterior domain $\Omega$, satisfying the uniform
$\delta-$cone condition, and a suitable subsequence of $\ep's$ (not
relabeled) such that
\bFormula{prop1}
| \Ome \setminus \Omega | + |\Omega \setminus \Ome| \to 0 \
\mbox{as}\ \ep \to 0.
\eF
For each $x_0 \in \partial \Omega$, there is $x_{\ep,0} \in \partial
\Ome$ such that $x_{\ep,0} \to x_0$, in particular,
\bFormula{prop2}
R^3 \setminus \Omega \subset B_d.
\eF
For any compact $K \subset \Omega$, there exists $\ep(K)$ such that
\bFormula{prop3}
K \subset \Ome \ \mbox{for all}\ \ep < \ep(K).
\eF

\end{itemize}

Property (\ref{prop1}) is important when studying stability of the
spectral properties of the Neumann Laplacian $\Den$ with respect to
$\ep$, see Arrieta and Krej\v ci\v r\' \i k \cite{ArrKre}. Note that
the limit domain \emph{need not} be of class $C^2$ but merely
Lipschitz, see Henrot and Pierre \cite[Theorem 2.4.7]{hepi05}.

\subsection{Main result}

Before stating our main result, we introduce the limit problem - the \emph{incompressible} Navier-Stokes system - satisfied by the limit velocity field
$\vc{U}$:

\bFormula{ins1}
\Div \vc{U} = 0,
\eF
\bFormula{ins2}
\Ov{\vr} \left( \partial_t \vc{U} + \Div (\vc{U} \otimes \vc{U} ) \right) + \Grad \Pi
= \mu \Delta \vc{U},
\eF
in $(0,T) \times \Omega$, with the condition at infinity
\bFormula{ins3}
|\vc{U}| \to 0 \ \mbox{as} \ |x| \to \infty,
\eF
and the initial condition
\bFormula{ins4}
\vc{U}(0, \cdot) = \vc{U}_0.
\eF

In the \emph{weak} formulation, the decay condition (\ref{ins3}) is
replaced by a single stipulation $\vc{U} \in L^2(0,T;
W^{1,2}(\Omega, R^3))$, the incompressibility constraint
(\ref{ins1}) is satisfied a.a.~in $(0,T) \times \Omega$, while
the momentum equation (\ref{ins2}), together with (\ref{ins4}), are
replaced by a family of integral identities;
\bFormula{ins5}
\int_0^T \intO{ \left( \Ov{\vr} \vc{U} \cdot \partial_t \varphi +
\Ov{\vr} \vc{U} \otimes \vc{U} : \Grad \varphi \right) } \ \dt
\eF
\[
 = \mu \int_0^T \intO{ \Grad \vc{U}:
\Grad \varphi } \ \dt - \intO{ \Ov{\vr} \vc{U}_0 \cdot \varphi (0, \cdot) },
\]
for any test function $\varphi \in \DC([0,T) \times \Omega; R^3)$ satisfying
$\Div \varphi = 0$.

Note that we have deliberately omitted to specify any boundary
conditions on $\partial \Omega$. In the low Mach number limit, we
can easily show that the impermeability condition $\vue \cdot
\vc{n}|_{\partial \Ome} = 0$ gives rise to the same property $\vc{U}
\cdot \vc{n}|_{\partial \Omega} = 0$ for the limit velocity field
whereas the boundary behavior of the tangential component of
$\vc{U}$ may be quite complex depending sensitively on the
asymptotic shape of the boundaries $\partial \Ome$, cf.
\cite{BuFeNe}. Sufficient conditions for $\vc{U}$ to satisfy the
no-slip condition
\bFormula{noslip}
\vc{U}|_{\partial \Omega} = 0
\eF
will be discussed in Section \ref{coa}.

Our main result reads as follows.

\bTheorem{m1}
Suppose that a family of domains $\{ \Ome \}_{\ep > 0}$ belongs to the class specified in Section \ref{dc}, with
\[
0 < \beta < \frac{1}{4},
\]
where $\ep^\beta$ is the radius of the balls in (\ref{dc1}). Let
$\{ \vre, \vue \}_{\ep > 0}$ be a family
weak solutions to the compressible Navier-Stokes system (\ref{i1} - \ref{i4a}) in $(0,T) \times \Ome$, supplemented with the initial conditions (\ref{m5}),
where
\[
r_{0,\ep} \equiv \frac{\vr_{0,\ep} - \Ov{\vr}}{\ep} \to r_0 \ \mbox{weakly in}\
L^2(R^3), \ \| r_{0,\ep} \|_{L^\infty(R^3)} \leq c, \ \Ov{\vr} > 0,
\]
\[
\vu_{0,\ep} \to \vc{u}_0 \ \mbox{weakly in} \ L^2(R^3; R^3),
\]
and the pressure satisfies (\ref{m1}) with $\gamma > 3/2$.

Then
\bFormula{mm1}
{\rm ess} \sup_{t \in (0,T)} \| \vre (t, \cdot) - \Ov{\vr} \|_{(L^2
+ L^q)(\Ome)} \to 0 \ \mbox{as} \ \ep \to 0 \ \mbox{for}\  1 \leq q
< \min \{ \gamma, 2 \},
\eF
\bFormula{mm2}
\| \vue \|_{L^2(0,T;W^{1,2}(\Ome;R^3))} \leq c,
\eF
and, at least for a suitable subsequence,
\bFormula{mm3}
\vue \to \vc{U} \
\mbox{in} \ L^2((0,T) \times K; R^3) \ \mbox{for any compact}\ K \subset \Ome,
\eF
where
\[
\vc{U} \in L^2(0,T; W^{1,2}(\Omega; R^3)), \ \vc{U} \cdot \vc{n}|_{\partial \Omega} = 0,
\]
is a weak solution of the incompressible Navier-Stokes system (\ref{ins1} - \ref{ins3}) in $(0,T) \times \Omega$, emanating from the initial data
\[
\vc{U}_0 = \vc{H} [\vc{u}_0],
\]
where $\Omega$ is the limit domain identified through (\ref{prop1} - \ref{prop3}) and $\vc{H}$ denotes the standard Helmholtz projection in $\Omega$.
\eT

The rest of the paper is devoted to the proof of Theorem \ref{Tm1}.
Note that  \emph{existence} of solution to the compressible Navier-Stokes system,
in the framework of weak solutions, was first established in the
seminal work of Lions \cite{LI4} for $\gamma \geq 9/5$, and the
result then extended in \cite{FNP} to the ``technically'' optimal
range $\gamma > 3/2$. A detailed discussion of various choices of
boundary conditions, including the case of exterior domains, may be
also found in the monograph by Novotn\' y and Stra\v skraba
\cite{NOST4}. The problem of limit boundary conditions, here
discussed in Section \ref{co}, was studied by Casado-Diaz,
Fernandez-Cara, and Simon \cite{CFS} in the periodic setting, and
later extended in \cite{BuFeNe}, \cite{BFNW}. Finally, we remark
there is an analogue of Theorem \ref{Tm1} for a family of
\emph{bounded} domains obtained by completely different methods (see
\cite{BuFe1}), based on  analysis of a viscous boundary layer
similar to Desjardins et al. \cite{DGLM}.

\section{Uniform bounds}
\label{ub}

All uniform bounds presented below may be viewed as a direct consequence of the energy inequality (\ref{m9}). To begin, similarly to \cite{FEINOV}, we introduce the \emph{essential} and \emph{residual} part of a function $h_\ep$ as
\[
[h_\ep]_{\rm ess} = \chi(\vre) h_\ep, \ [h_\ep]_{\rm res} = h - h_{\rm ess},
\]
where
\[
\chi \in \DC(0, \infty), \ 0 \leq \chi \leq 1, \ \chi \equiv 1 \ \mbox{in an open neighborhood of} \ \Ov{\vr}.
\]

\subsection{Energy bounds}

Since the initial data satisfy (\ref{m5}), (\ref{m6}),
the initial energy $E_\ep (\vr_{0,\ep}, \vu_{0,\ep})$ in (\ref{m9}) remains bounded
uniformly for $\ep \to 0$, where we have used hypothesis (\ref{m1}).
Consequently, we deduce the following list of estimates:
\bFormula{ub4}
{\rm ess} \sup_{t \in (0,T)} \intOe{ \vre |\vue|^2 } \leq c,
\eF
\bFormula{ub1}
{\rm ess} \sup_{t \in (0,T)} \intOe{ \left[ \frac{\vre - \Ov{\vr}}{\ep} \right]_{\rm ess}^2 } \leq c,\
{\rm ess} \sup_{t \in (0,T)} \intOe{ \left[ \vre \right]_{\rm res}^\gamma } \leq \ep^2 c,
\
{\rm ess} \sup_{t \in (0,T)} \intOe{  1_{\rm res} } \leq \ep^2 c,
\eF
and
\bFormula{ub5}
\int_0^T \intOe{ \left| \Grad \vue + \Grad^t \vue - \frac{2}{3} \Div \vue \tn{I} \right|^2
} \ \dt \leq c,
\eF
where the constants are independent of $\ep \to 0$.

Moreover, relation (\ref{ub1}) immediately yields
\bFormula{ae6}
{\rm ess} \sup_{t \in (0,T)} \left\| \left[ \frac{\vre - \Ov{\vr}
}{\ep} \right]_{\rm res} \right\|_{L^q(\Ome)} \leq c
\ep^{\frac{2-q}{q}} \ \mbox{for any}\ 1 \leq q \leq \min\{ \gamma, 2
\},
\eF
which, together with (\ref{ub1}) gives rise to (\ref{mm1}).

Finally, since the family of domains $\{ \Ome \}_{\ep > 0}$ admits
the uniform Korn's inequality (\ref{uki}), we can combine (\ref{ub1}), (\ref{ub5})
to conclude that
\bFormula{ub6}
\| \vue \|_{L^2(0,T; W^{1,2}(\Ome; R^3))} \leq c,
\eF
uniformly for $\ep \to 0$.

\subsection{Convergence}

As the family $\{ \Ome \}_{\ep > 0}$ possesses the uniform extension property
(\ref{uep}), we may assume that
\bFormula{convv}
\vue \to \vc{U} \ \mbox{weakly in} \ L^2(0,T; W^{1,2}(R^3; R^3)).
\eF

Moreover, by virtue of the uniform bounds (\ref{ub1}), (\ref{ub6}),
we can perform the limit in (\ref{m7}) to obtain
\[
\Div \vc{U} = 0 \ \mbox{a.a. in} \ (0,T) \times \Omega,
\]
and, similarly, one can pass to the limit in the weak formulation of
momentum equation (\ref{m8}) to deduce that
\[
\int_0^T \intO{ \left( \Ov{\vr} \vc{U} \cdot \partial_t \varphi +
\Ov{\vr \vc{u} \otimes \vc{u}} : \Grad \varphi \right) } \ \dt
\]
\[
 = \mu \int_0^T \intO{ \Grad \vc{U}:
\Grad \varphi } \ \dt - \intO{ \Ov{\vr} \vc{U}_0 \cdot \varphi (0, \cdot) },
\]
for any test function $\varphi \in \DC([0,T) \times \Omega; R^3)$ satisfying
$\Div \varphi = 0$, where the symbol $\Ov{\vr \vu \otimes \vu}$ denotes a weak limit
of $\{ \vre \vue \otimes \vue \}_{\ep > 0}$. Accordingly, in order to finish the proof of Theorem \ref{Tm1}, we have to show
\[
\Ov{ \vr \vu \otimes \vu } = \Ov{\vr} \vc{U} \otimes \vc{U},
\]
or, equivalently, the strong convergence of the velocities claimed in (\ref{mm3}).
This will be our goal in the remaining part of the paper.
We remark that, by virtue of (\ref{prop1} - \ref{prop3}), it is easy to check that
the the limit velocity field satisfies the impermeability condition
\[
\vc{U} \cdot \vc{n}|_{\partial \Omega} = 0,
\]
in a weak sense.

\section{The acoustic equation}
\label{ae}

In this section, we introduce a weak formulation of the acoustic equation (\ref{i11} - \ref{i13}) and discuss its basic properties.

\subsection{Weak formulation}

With
\[
r_\ep \equiv \frac{\vre - \Ov{\vr}}{\ep}, \ \vc{V}_\ep \equiv \vre \vue,
\]
the Navier-Stokes system (\ref{m7}), (\ref{m8})
can be written in the form:
\bFormula{ae1}
\int_0^T \intOe{ \Big( \ep r_\ep \partial_t \varphi + \vc{V}_\ep \cdot \Grad \varphi \Big) } \ \dt = - \intOe{ \ep r_{0,\ep} \varphi (0, \cdot) },
\eF
for any $\varphi \in \DC([0,T) \times \Ov{\Omega})$,
\bFormula{ae2}
\int_0^T \intOe{ \Big( \ep \vc{V}_\ep \cdot \partial_t \varphi + p'(\Ov{\vr} ) r_\ep \Div
\varphi \Big) } \ \dt  = - \intOe{ \ep \vr_{0,\ep} \vu_{0,\ep} \cdot \varphi(0, \cdot)}
\eF
\[
+ \int_0^T \intOe{ \Big( \ep \tn{S} (\Grad \vue) : \Grad \varphi -
\ep \vre \vue \otimes \vue : \Grad \varphi \Big)} \ \dt
\]
\[
- \int_0^T \intOe{ \ep \left(  \frac{1}{\ep^2} \Big( p (\vre) -
p'(\Ov{\vr}) (\vre - \Ov{\vr}) - p(\Ov{\vr}) \Big) \Div \varphi
\right)   } \ \dt,
\]
for any $\varphi \in \DC([0,T) \times \Ov{\Omega}_\ep;R^3)$, $\varphi \cdot \vc{n}|_{\partial \Ome} = 0$.

Now, thanks to the slip boundary condition (\ref{i4}), observe that
$\Grad \Delta^{-1}_{\ep, N}[\varphi]$, where
$\Delta_{\ep, N}$ is the Neumann Laplacian in $\Ome$, is an
admissible test function in (\ref{ae2}). Consequently, we obtain
\bFormula{ae3}
\int_0^T \intOe{ \Big( \ep \Phi_\ep  \partial_t \varphi -
p'(\Ov{\vr} ) r_\ep \varphi \Big) } \ \dt  =  - \intOe{ \ep
\vr_{0,\ep} \vu_{0,\ep} \cdot \Grad \Delta^{-1}_{\ep,N} [\varphi(0,
\cdot)]},
\eF
\[
- \int_0^T \intOe{ \Big( \ep \tn{S} (\Grad \vue) : \Grad^2
\Delta^{-1}_{\ep,N} [\varphi] - \ep \vre \vue \otimes \vue : \Grad^2
\Delta^{-1}_{\ep, N}[\varphi] \Big)} \ \dt
\]
\[
+ \int_0^T \intOe{ \ep \left(  \frac{1}{\ep^2} \Big( p (\vr) -
p'(\Ov{\vr}) (\vr - \Ov{\vr}) - p(\Ov{\vr}) \Big) \varphi \right) }
\ \dt ,
\]
for all $\varphi \in \DC([0,T) \times \Ov{\Omega})$, where $\Phi_\ep$ is the acoustic potential, meaning,
\[
\vc{V}_\ep = \vc{H}_\ep [\vc{V}_\ep] + \Grad \Phi_\ep,
\]
where $\vc{H}_\ep$ denotes the standard Helmholtz projection in
$\Omega_\ep$. Note that for $\varphi \in \DC([0,T) \times
\Ov{\Omega}_\ep)$ and $\varepsilon > 0$ fixed, the test function
$\Grad \Delta^{-1}_{\ep,N}[\varphi]$ is continuously differentiable
in $\Ov{\Omega}_\ep$ and belongs to the space $\DC([0,T);
W^{1,2}(\Ome;R^3))$.

The Helmholtz projection $\vc{H}_\ep[ \vc{v}]$ of a function $\vc{v} \in L^p(\Ome; R^3)$ is defined as
\[
\vc{H}_\ep [\vc{v}] = \vc{v} - \Grad \Phi,
\]
where $\Phi \in D^{1,p} ({\Omega}_\ep)$ is the unique solution of the problem
\[
\intOe{ \Grad \Phi \cdot \Grad \varphi } = \intOe{ \vc{v} \cdot
\Grad \varphi } \ \mbox{for all} \ \varphi \in  \DC (\Ov{\Omega}_\ep
),
\]
in other words, at least formally,
\[
\Delta \Phi = \Div \vc{v} \ \mbox{in} \in \Ome, \ \Grad \Phi \cdot
\vc{n}|_{\partial \Ome} = \vc{v} \cdot \vc{n}|_{\partial \Ome} , \
|\Phi| \to 0 \ \mbox{for}\ |x| \to \infty.
\]
Here, the symbol $D^{1,p}(\Ome)$ denotes the homogeneous Sobolev
space - a completion of $\DC(\Ov{\Omega}_\ep )$ with respect to the
norm $\| \Grad \Phi \|_{L^p(\Ome;R^3)}$. We have Sobolev's
inequality
\bFormula{si}
\| \Phi \|_{L^q(\Ome)} \leq c(p) \| \Grad \Phi \|_{L^p(\Ome;R^3)},\ q = \frac{3p}{3 - p}
\ \mbox{for any} \ 1 \leq p < 3,
\eF
for any $\Phi \in D^{1,p}(\Ome)$,
where, since $\{ \Ome \}_{\ep > 0}$ admits the uniform extension property (\ref{uep}), the
constant $c(p)$ is independent of $\ep$.

Finally, the equation (\ref{ae1}) reads
\bFormula{ae4}
\int_0^T \intOe{ \Big( \ep r_\ep \partial_t \varphi + \Grad \Phi_\ep \cdot \Grad \varphi \Big) } \ \dt = - \intOe{ \ep r_{0,\ep} \varphi (0, \cdot) }
\eF
for any $\varphi \in \DC([0,T) \times \Ov{\Omega})$. Equations (\ref{ae3}), (\ref{ae4}) represent a weak formulation of the acoustic equation (\ref{i11}), (\ref{i12}), with the Neumann boundary condition (\ref{i13}) implicitly included through the class of test functions in (\ref{ae4}).

\subsection{Uniform bounds, part I}

It follows from the uniform bounds established in
(\ref{ub4}), (\ref{ub1}), (\ref{ae6}) that
\[
r_\ep = [r_\ep ]_{\rm ess} + [r_\ep ]_{\rm res},
\]
satisfies
\bFormula{ae5}
{\rm ess} \sup_{t \in (0,T)} \left\| \left[ r_\ep \right]_{\rm ess} \right\|_{L^2(\Ome)} \leq c,
\eF
and
\bFormula{ae5a}
{\rm ess} \sup_{t \in (0,T)} \left\| \left[ r_\ep \right]_{\rm res}
\right\|_{L^q(\Ome)} \leq \ep^{ \frac{2 - q}{ q}} c \ \mbox{for any}
\ 1 \leq q < \min\{ \gamma, 2 \}.
\eF

Similarly,
\[
\vc{V}_\ep = [ \vc{V}_\ep ]_{\rm ess} + [ \vc{V}_\ep ]_{\rm res},
\]
where, in accordance with (\ref{ub4}), (\ref{ub1}),
\bFormula{ae7}
{\rm ess} \sup_{t \in (0,T)} \|[ {\vc{V}_\ep } ]_{\rm ess} \|_{L^2(\Ome;R^3)}
=
{\rm ess} \sup_{t \in (0,T)} \|[ \sqrt{\vre} ]_{\rm ess} \sqrt{\vre}
\vue \|_{L^2(\Ome;R^3)} \leq c,
\eF
and
\bFormula{ae8}
{\rm ess} \sup_{t \in (0,T)} \|[ {\vc{V}_\ep } ]_{\rm res} \|_{L^q(\Ome;R^3)}
=
{\rm ess} \sup_{t \in (0,T)} \|[ \sqrt{\vre} ]_{\rm res} \sqrt{\vre}
\vue \|_{L^q(\Ome;R^3)} \leq c \ep^{1/\gamma}, \ q =
\frac{2\gamma}{\gamma + 1}.
\eF

It remains to find suitable bounds on the forcing terms in acoustic
equation (\ref{ae3} - \ref{ae4}). To this end, we need the elliptic
estimates for the Neumann Laplacian $\Den$ discussed in the next
section.

\subsection{Elliptic estimates and Helmholtz decomposition in $\Ome$}

In order to control the forcing terms as well as the initial data in
the acoustic equation, we need bounds on $\Grad^2 v$ in terms of $\Den [v]$. As the curvature,
represented by the radius of the balls in (\ref{dc1}),
\emph{is not} uniformly bounded, the $W^{2,p}-$elliptic bounds may ``blow-up'' for $\ep \to 0$.

\subsubsection{$W^{2,p}-$bounds}

In order to obtain $W^{2,p}-$bounds, we consider the rescaled family of domains
\bFormula{sd}
\widehat \Ome \equiv \frac{1}{\ep^\beta} \Ome,
\eF
where the exponent $\beta > 0$ is the same as in (\ref{dc1}). Accordingly, the rescaled
domains $\widehat \Ome$ are of \emph{uniform} $C^2$-class, in particular, the standard
elliptic theory yields
\bFormula{ee1}
\| \Grad^2 v \|_{L^p(\widehat{\Omega}_\ep; R^{3 \times 3})} \leq c(p) \left( \| \Delta_x v \|_{L^p(\widehat{\Omega}_\ep )} +
\| v \|_{L^p(\widehat{\Omega}_\ep)} \right) \ \mbox{for any}\ 1 < p < \infty,
\eF
for any $v \in \DC(\Ov{\widehat{\Omega}}_\ep )$ satisfying $\Grad v
\cdot \vc{n}|_{\partial \widehat{\Omega}_\ep } = 0$. It is important
to notice that, by virtue of the hypotheses introduced in Section
\ref{dc}, the constant $c(p)$ depends only on the rescaled radius
$c_b$ of the balls appearing in (\ref{dc1}).

Consequently, returning to the original domains $\Ome$ we may infer that
\bFormula{ee2}
\| \Grad^2 v \|_{L^p(\Ome; R^{3 \times 3})} \leq c(p) \left( \| \Delta_x v \|_{L^p(\Ome)} + \frac{1}{\ep^{2 \beta}}
\| v \|_{L^p(\Ome)} \right) \ \mbox{for}\ 1 < p < \infty,
\eF
for any $v \in \DC(\Ov{\Omega}_\ep )$ such that $\Grad v \cdot
\vc{n}|_{\partial {\Omega}_\ep } = 0$, with $c(p)$ independent of
$\ep \to 0$.

\subsubsection{Helmholtz decomposition}

Consider the family of rescaled domains $\widehat{\Omega}_\ep$ introduced in (\ref{sd}),
with the associated Helmholtz projections $\widehat{\vc{H}}_\ep$. By virtue of the
result by Farwig, Kozono, and Sohr \cite{FAKOSO}, we have
\bFormula{ee3}
\| \widehat{\vc{H}}_\ep [\vc{v}] \|_{(L^p \cap L^2) (\widehat{\Omega}_\ep, R^3)} \leq
c(p) \| \vc{v} \|_{(L^p \cap L^2) (\widehat{\Omega}_\ep, R^3)} \ \mbox{for any}\
2 \leq p < \infty,
\eF
where, similarly to the preceding part, the constant $c(p)$ depends
only on $c_b$. Going back to the original domain $\Ome$ we therefore
obtain
\bFormula{ee4}
\| {\vc{H}}_\ep [\vc{v}] \|_{(L^p \cap L^2) ({\Omega}_\ep, R^3)} \leq
\ep^{-\beta\left( \frac{3}{2} - \frac{3}{p} \right)} c(p) \| \vc{v} \|_{(L^p \cap L^2) ({\Omega}_\ep, R^3)} \ \mbox{for any}\
2 \leq p < \infty,
\eF
uniformly for $\ep \to 0$.

Similarly, by means of a duality argument,
\bFormula{ee5}
\| {\vc{H}}_\ep [\vc{v}] \|_{(L^p + L^2) ({\Omega}_\ep, R^3)} \leq
\ep^{-\beta\left( \frac{3}{p} - \frac{3}{2} \right)} c(p) \| \vc{v} \|_{(L^p + L^2) ({\Omega}_\ep, R^3)} \ \mbox{for any}\
1 < p < 2.
\eF

The estimates (\ref{ee4}), (\ref{ee5}) arise also in problems of
homogenization  and need not be optimal, cf. Masmoudi \cite{MAS3},
\cite{MASCras}. On the other hand, as the limit domain
$\Omega$ may be only Lipschitz, it is not surprising that the
$L^p-$bounds in (\ref{ee4}), (\ref{ee5}) blow up for $\ep \to 0$.

\subsection{Uniform bounds, part II}

With the bounds established in the previous section at hand,
we are able to control the forcing terms  in the acoustic equation (\ref{ae3}).
To begin, relation (\ref{ee2}) implies that
\[
\left| \intOe{  \tn{S} (\Grad \vue) : \Grad^2 \Delta^{-1}_{\ep,N}
[\varphi] } \right|
\]
\[
\leq  c \| \tn{S} (\Grad \vue )
\|_{L^2(\Ome;R^{3 \times 3})} \left( \| \varphi \|_{L^2(\Ome)} + \frac{1}{\ep^{2 \beta}} \| (-\Den)^{-1}[\varphi] \|_{L^2(\Ome)}
\right),
\]
therefore, by means of the Riesz representation theorem,
\bFormula{ae11}
\int_0^T \intOe{   \tn{S} (\Grad \vue) : \Grad^2 \Delta^{-1}_{\ep,N}
[\varphi] } \ \dt = \frac{1}{\ep^{ 2 \beta}} \int_0^T \intOe{ \left( F^1_\ep \varphi
+ F^2_\ep (-\Den)^{-1} [\varphi ] \right) }
\ \dt ,
\eF
where
\bFormula{ae12}
\| F^i_\ep \|_{L^2((0,T) \times \Ome)} \leq c,\ i=1,2, \ \mbox{uniformly
for}\ \ep \to 0.
\eF

Similarly, we can write
\[
\intOe{ \vre \vue \otimes \vue : \Grad^2 \Delta^{-1}_{\ep,
N}[\varphi] }
\]
\[
= \intOe{ [\vre]_{\rm ess} \vue \otimes \vue : \Grad^2
\Delta^{-1}_{\ep, N}[\varphi] } + \intOe{ [\sqrt{\vre}]_{\rm res}
\sqrt{\vre} \vue \otimes \vue : \Grad^2 \Delta^{-1}_{\ep,
N}[\varphi] },
\]
where, by virtue of (\ref{ee2}),
\[
\left| \intOe{ [\vre]_{\rm ess} \vue \otimes \vue : \Grad^2
\Delta^{-1}_{\ep, N}[\varphi] } \right|
\]
\[
\leq \| [\vre]_{\rm ess} \vue \|_{L^2(\Ome, R^3)} \| \vue
\|_{L^6(\Ome;R^3)} \| \Grad^2 \Delta^{-1}_{\ep, N}[\varphi]
\|_{L^3(\Ome;R^{3 \times 3})}
\]
\[
\leq c \| [\vre]_{\rm ess} \vue \|_{L^2(\Ome, R^3)} \| \vue
\|_{L^6(\Ome;R^3)} \left( \| \varphi \|_{L^3(\Ome)} + \frac{1}{\ep^{2 \beta} }
\| (-\Den)^{-1} [\varphi ] \|_{L^3(\Ome)} \right),
\]
and, by interpolation, the uniform extension property, and Sobolev's inequality,
\[
\| \varphi \|_{L^3(\Ome)} \leq c_1 \left( \| \varphi \|_{L^2(\Ome)}  + \| \varphi \|_{L^6(\Ome)} \right)
\]
\[
\leq c_2 \left( \| \varphi \|_{L^2(\Ome)}  + \| \Grad \varphi \|_{L^2(\Ome)} \right) =
c_2 \left( \| \varphi \|_{L^2(\Ome)}  + \| (-\Den)^{1/2}[ \varphi ] \|_{L^2(\Ome)} \right),
\]
while, by the same token,
\[
\| (-\Den )^{-1} [\varphi] \|_{L^3(\Ome)} \leq c_1 \left( \| (-\Den)^{-1}[\varphi] \|_{L^2(\Ome)}  + \| (-\Den)^{-1/2}[\varphi] \|_{L^2(\Ome)} \right).
\]
Consequently, there exist functions $F^i_\ep$, $i=3, \dots, 6$,
\bFormula{ae13}
\| F^i_\ep \|_{L^2((0,T) \times \Ome)} \leq c,\ i=3, \dots, 6 \
\mbox{uniformly for}\ \ep \to 0,
\eF
such that
\bFormula{ae14}
\int_0^T \intOe{ [\vre]_{\rm ess} \vue \otimes \vue : \Grad^2
\Delta^{-1}_{\ep, N}[\varphi] } \ \dt
\eF
\[
= \frac{1}{\ep^{2\beta}} \int_0^T \intOe{ \left( F^3_\ep \varphi +
F^4_\ep (-\Den)^{-1/2}[\varphi] + F^5_\ep (-\Den)^{1/2}[\varphi] +
F^6_\ep (-\Den)^{-1}[\varphi] \right) } \ \dt.
\]

Furthermore, in accordance with (\ref{ae8}),
\[
\left| \intOe{ [\vre]_{\rm res} \vue \otimes \vue : \Grad^2
\Delta^{-1}_{\ep, N}[\varphi] } \right| \leq \ep^{1/\gamma} c \| \vu
\|_{L^6(\Ome;R^3)} \| \Grad^2 \Delta^{-1}_{\ep, N}[\varphi]
\|_{L^r(\Ome; R^{3 \times 3})},
\]
where, by virtue of (\ref{ee2}),
\[
\| \Grad^2 \Delta^{-1}_{\ep, N}[\varphi]
\|_{L^r(\Ome; R^{3 \times 3})} \leq c \left( \| \varphi \|_{L^r(\Ome)} + \frac{1}{\ep^{2 \beta}} \| (-\Den)^{-1}[ \varphi ] \|_{L^r(\Ome)} \right),
\]
with
\[
\frac{\gamma + 1}{2 \gamma} + \frac{1}{6} + \frac{1}{r} = 1.
\]

In what follows, we suppose $\gamma < 2$, $r \geq 3$ as, otherwise,
the estimates are the same as in (\ref{ae14}). Thus, applying once
more (\ref{ee2}), we obtain
\[
\| \varphi \|_{L^r(\Ome)} \leq c \left( \| (-\Den) [ \varphi ] \|_{L^2(\Ome)} +
\frac{1}{\ep^{2 \beta}} \| \varphi \|_{L^2(\Ome)} \right),
\]
and, similarly,
\[
\| (-\Den)^{-1} \varphi \|_{L^r(\Ome)} \leq c \left( \|  \varphi  \|_{L^2(\Ome)} +
\frac{1}{\ep^{2 \beta}} \| (-\Den)^{-1} [\varphi] \|_{L^2(\Ome)} \right).
\]
Consequently, we get the same result as (\ref{ae14}) provided $\max
\{ \frac{1}{\gamma} , \frac{1}{2} \} \geq 2 \beta$, in particular if
$\beta < 1/4$.

Thus, we may infer that
\bFormula{ae15}
\int_0^T \intOe{ \vre \vue \otimes \vue : \Grad^2
\Delta^{-1}_{\ep, N}[\varphi] } \ \dt
\eF
\[
= \frac{1}{\ep^{2\beta}} \int_0^T \int_{\Omega_\ep} \left( F^3_\ep
\varphi + F^4_\ep (-\Den)^{-1/2}[\varphi] + F^5_\ep
(-\Den)^{1/2}[\varphi] \right.
\]
\[
\left. + F^6_\ep (-\Den)[\varphi] + F^7_\ep (-\Den)^{-1}[\varphi]
\right) \dx \ \dt,
\]
where $F^i_\ep$ satisfy (\ref{ae13}). Note that the same symbol
$F^i_\ep$ may stand for different functions than above.

Finally,
\[
\left| \intOe{ \left(  \frac{1}{\ep^2} \Big( p (\vre) - p'(\Ov{\vr})
(\vre - \Ov{\vr}) - p(\Ov{\vr}) \Big) \varphi \right) } \right|
\]
\[
\leq \left\| \frac{1}{\ep^2} \Big( p (\vre) - p'(\Ov{\vr}) (\vre -
\Ov{\vr}) - p(\Ov{\vr}) \Big) \right\|_{L^1(\Ome)} \| \varphi
\|_{L^\infty(\Ome)},
\]
where, in accordance with the uniform extension property,
\[
\| \varphi \|_{L^\infty(\Ome)} \leq c \left( \| \Grad \varphi
\|_{L^6(\Ome, R^{3})} + \| \varphi \|_{L^6(\Ome)} \right) \leq c
\left( \| \Grad^2 \varphi \|_{L^2(\Ome, R^{3 \times 3})} + \| \Grad
\varphi \|_{L^2(\Ome;R^3)} \right),
\]
and, by virtue of (\ref{ee2}),
\[
\| \Grad^2 \varphi \|_{L^2(\Ome, R^{3 \times 3})} \leq
c \left( \| (-\Den) [\varphi] \|_{L^2(\Ome)} + \frac{1}{\ep^{2 \beta}} \| \varphi \|_{L^2(\Ome)}
\right),
\]
\[
\| \Grad \varphi \|_{L^2(\Ome;R^3)} = \| (-\Den)^{1/2} \varphi \|_{L^2(\Ome)}.
\]

Seeing that, by virtue of the uniform bounds established (\ref{ub1}), (\ref{ae6}),
\[
{\rm ess} \sup_{t \in (0,T)} \left\| \frac{1}{\ep^2} \Big( p (\vre) - p'(\Ov{\vr}) (\vre - \Ov{\vr}) - p(\Ov{\vr}) \Big) \right\|_{L^1(\Ome)} \leq c,
\]
we conclude that
\bFormula{ae18}
\int_0^T
\intOe{ \left(  \frac{1}{\ep^2} \Big( p (\vre) - p'(\Ov{\vr})
(\vre - \Ov{\vr}) - p(\Ov{\vr}) \Big) \varphi \right) } \ \dt
\eF
\[
= \frac{1}{\ep^{2\beta}} \int_0^T \intOe{ \left( F^8_\ep \varphi +
F^9_\ep (-\Den)^{1/2}[\varphi] + F^{10}_\ep (-\Den)[\varphi] \right)
} \ \dt,
\]
where
\bFormula{ae18a}
\| F^i_\ep \|_{L^2((0,T) \times \Ome)} \leq c,\ i=8, 9, 10 \
\mbox{uniformly for}\ \ep \to 0.
\eF

\subsection{The acoustic equation revisited}
Using relations (\ref{ae7}), (\ref{ae8}), together with (\ref{ee5}),
we can  write the acoustic potential $\Phi_\ep$ in
the form
\[
\Phi_\ep = \Phi^1_\ep + \Phi^2_\ep,
\]
where
\bFormula{ae9}
{\rm ess} \sup_{t \in (0,T)} \| \Phi^1_\ep \|_{D^{1,2}(\Ome; R^3)}
\leq c,
\eF
\bFormula{ae10}
{\rm ess} \sup_{t \in (0,T)} \| \Phi^2_\ep \|_{L^2(\Ome; R^3)} \to 0
\ \mbox{as}\ \ep \to 0.
\eF

In view of the uniform bounds obtained in the previous section,
the acoustic equation (\ref{ae3}), (\ref{ae4}) can be written in the
concise form
\bFormula{ae19}
\int_0^T \intOe{ \Big( \ep r_\ep \partial_t \varphi + \Grad \Phi_\ep
\cdot \Grad \varphi \Big) } \ \dt = - \intOe{ \ep r_{0,\ep} \varphi
(0, \cdot) }
\eF
for any $\varphi \in \DC([0,T) \times \Ov{\Omega})$, and
\bFormula{ae20}
\int_0^T \intOe{ \Big( \ep \Phi_\ep  \partial_t \varphi -
p'(\Ov{\vr} ) r_\ep \varphi \Big) } \ \dt  =  - \intOe{ \ep \Phi_{0,\ep}
\varphi(0, \cdot)}
\eF
\[
+ \ep^{1 - 2 \beta} \int_0^T \int_{\Omega_\ep} \Big( G^1_\ep \varphi + G^2_\ep
(-\Den)^{-1/2} [ \varphi] + G^3_\ep (- \Den)^{1/2} [\varphi]
\]
\[
+ G^4_\ep (-\Den)[\varphi] +G^5_\ep (-\Den)^{-1}[\varphi] \Big) \dx \ \dt,
\]
for all $\varphi \in \DC([0,T) \times \Ov{\Omega})$, $\Grad \varphi
\cdot \vc{n}|_{\partial \Ome} = 0$, where
\bFormula{ae20a}
\| r_{0, \ep} \|_{L^2(\Ome)} + \| (- \Den )^{-1/2}[\Phi_{0, \ep}] \|_{L^2(\Ome)}
\leq c,
\eF
\bFormula{ae21}
\| G^i_\ep \|_{L^2(\Ome)} \leq c,\ i=1, \dots, 5,
\eF
uniformly for $\ep \to 0$.

\section{Spectral analysis of Neumann Laplacian on varying domains}
\label{nl}

As  observed, the Neumann Laplacian $\Den$ plays a crucial
role in the analysis of acoustic waves. We recall that $-\Den$ may
be viewed as a non-negative self-adjoint operator on the space
$L^2(\Ome)$, with
\[
{\cal D}(- \Den) = \left\{ w \in W^{1,2}(\Ome) \ \Big| \ \intOe{
\Grad w \cdot \Grad \varphi } = \intOe{ g \varphi } \right.
\]
\[
\left. \mbox{for all}\ \varphi \in \DC(\Ov{\Omega}_\ep)  \ \mbox{and
a certain}\ g \in L^2(\Ome) \right\},\ -\Den w = g.
\]
Since the boundaries $\partial \Ome$ are regular, the standard
elliptic theory yields
\[
{\cal D}(-\Den) = \left\{ w \in W^{2,2}(\Ome) \ \Big| \ \Grad w
\cdot \vc{n}|_{\partial \Ome} = 0 \right\}.
\]
We denote by $\{ {\cal P}_{\ep, \lambda} \}_{\lambda \geq 0 }$ the
associated family of spectral projections. The following analysis is
a slight modification of \cite[Section 2]{EF84}, similar problems
were studied by Rauch and Taylor \cite{RauTay}.

\subsection{Spectral measures}

\label{sa}

Our goal is to obtain dispersive estimates, and, in
particular, local decay of acoustic waves, with a rate
\emph{independent} of the scaling parameter $\ep$. To this end, we
introduce the \emph{spectral measure} $\mu_{\ep, \varphi}$
associated to a function $\varphi \in L^2(\Ome)$ through Stone's
formula (see Reed and Simon \cite[Theorem VII.13]{ReSi1})
\bFormula{sa1}
\mu_{\ep, \varphi}(a,b)
\eF
\[
= \lim_{\delta \to 0+} \lim_{\eta \to 0+} \int_{a + \delta}^{b -
\delta} \left< \left( \frac{1}{ - \Den - \lambda - {\rm i} \eta} -
\frac{1}{- \Den - \lambda + {\rm i} \eta} \right)[\varphi]; \varphi
\right>_{\Ome} \ {\rm d}\lambda,
\]
where the symbol $\left< \cdot; \cdot \right>_{\Ome}$ denotes the
standard (complex) scalar product in $L^2(\Ome)$.

Now, the crucial observation is that we may perform the limit $\eta
\to 0+$ in (\ref{sa1}) as soon as $\varphi \in \DC(\Ov{\Omega}_\ep
)$ since the operators $-\Den$ satisfy the limiting absorption principle, see Leis
\cite{Leis}, Va{\u\i}nberg \cite[Chapter VIII]{Vain}:

\medskip

\centerline{\bf Limiting absorption principle (LAP)}

\medskip \noindent
{\it The operators $(1 + |x|^2)^{-s/2} \circ (- \Den - \lambda \pm
{\rm i \eta )^{-1} \circ (1 + |x|^2)^{-s/2}}$ are bounded on
$L^2(\Ome)$ for any $s
> 1$ uniformly for $\lambda$ belonging to compact subsets of
$(0,\infty)$ and $\eta > 0$.}

\medskip

Consequently, the spectral measure $\mu_{\ep, \varphi}$, $\varphi
\in \DC(\Ov{\Omega}_\ep)$, is absolutely continuous with respect to
the Lebesgue measure on $(0, \infty)$ and Stone's formula
(\ref{sa1}) reduces to
\bFormula{sa2}
\mu_{\ep,\varphi}(a,b) = \int_a^b \left< \left( w^-_{\ep, \lambda} -
w^+_{\ep, \lambda} \right) ;\varphi \right>_{\Ome} \ {\rm d}
\lambda,\ 0 < a < b,
\eF
where $w^{\pm}_{\ep, \lambda}$ are solutions of the Neumann problem
\bFormula{sa3}
\Delta w^{\pm}_{\ep, \lambda} + \lambda w^{\pm}_{\ep, \lambda} =
\varphi \ \mbox{in}\ \Ome,\ \Grad w^{\pm}_{\ep, \lambda} \cdot
\vc{n}|_{\partial \Ome} = 0,
\eF
uniquely determined by Sommerfeld's radiation condition
\bFormula{sa4}
\lim_{r \to \infty} r \left( \partial_r \pm {\rm i} \sqrt{\lambda}
\right) w^{\pm}_{\ep, \lambda} = 0,\ r \equiv |x|,
\eF
see Va{\u\i}nberg \cite[Chapter VIII]{Vain}.

Our goal is a uniform bound on the norm of the functions
$w^{\pm}_{\ep, \lambda}$ independent of the scaling parameter $\ep$.
To this end, we fix $R > 2d$ so that $\partial \Ome \subset B_{R/2}$
for all $\ep$. Following \cite{EF84} we recall the explicit formula
for the exterior Dirichlet problem
\[
\Delta v^{\pm}_{\ep, \lambda} + \lambda v^{\pm}_{\ep, \lambda} = 0 \
\mbox{in}\ R^3 \setminus \Ov{B}_R,\ v^{\pm}_{\ep, \lambda}
|_{\partial B_R} = \tilde v^{\pm}_{\ep, \lambda},
\]
supplemented with the radiation condition (\ref{sa4}), namely
\bFormula{sa5}
v^{\pm}_{\ep, \lambda} (x) = \sum_{l=0}^\infty \sum_{m=-l}^l a^m_l
Y^m_l (\theta, \phi) \frac{h^{(1)}_l (\pm \sqrt{\lambda}r)
}{h^{(1)}_l (\pm \sqrt{\lambda} R)} \ \mbox{for all}\ x \in R^3
\setminus \Ov{B}_{R},
\eF
where $(r,\theta,\phi)$ are polar coordinates, $Y^m_l$ are spherical
harmonics of order $l$, $h^{(1)}_l$ are spherical Bessel functions,
and
\[
\tilde v^{\pm}_{\ep, \lambda} (x) = \sum_{l=0}^\infty \sum_{m=-l}^l
a^m_l Y^m_l (\theta,\phi) \ \mbox{for}\ |x| = R,
\]
see Chandler-Wilde and Monk \cite{ChWMon}, N\' ed\' elec
\cite{Nedel}.

Assume that $\varphi \in \DC(\Ome)$ in (\ref{sa3}) is such that
\[
{\rm supp}[ \varphi ] \subset B_R.
\]
We claim that the functions $w^{\pm}_{\ep, \lambda}$ solving
(\ref{sa3}), (\ref{sa4}) admit a uniform bound
\bFormula{sa6}
\| w^{\pm}_{\ep, \lambda} \|_{L^2(B_{3R} \cap \Ome)} \leq c \|
\varphi \|_{L^2(\Ome)},
\eF
with $c$ independent of $\ep \to 0$, provided $\lambda$ belongs to a
compact subinterval of $(0, \infty)$. In order to see (\ref{sa6}),
we argue by contradiction assuming the existence of sequences
\[
\varphi_\ep \in \DC( B_R \cap \Omega_\ep),\ \| \varphi_\ep
\|_{L^2(B_R \cap \Omega_\ep)} = 1,
\]
\[
\lambda_\ep \to \lambda \in (0, \infty),
\]
such that the corresponding (unique) solutions $w^{\pm}_{\ep,
\lambda_\ep}$ of (\ref{sa3}), (\ref{sa4}) satisfy
\[
\| w^{\pm}_{\ep, \lambda_\ep} \|_{L^2(B_{3R} \cap \Omega_\ep)} \to
\infty \ \mbox{for}\  \ep \to 0.
\]

Setting
\[
v^{\pm}_{\ep, \lambda_\ep} = \frac{ w^{\pm}_{\ep, \lambda_\ep} } {\|
w^{\pm}_{\ep, \lambda_\ep} \|_{L^2(B_{3R} \cap \Omega_\ep)} },
\]
we check that $v^{\pm}_{\ep, \lambda_\ep}$ is the unique solution of
the problem
\[
\Delta v^{\pm}_{\ep,\lambda} + \lambda_\ep v^{\pm}_{\ep,\lambda} =
\frac{\varphi_\ep}{ \| w^{\pm}_{\ep, \lambda_\ep} \|_{L^2(B_{3R}
\cap \Omega_\ep)} } \ \mbox{in}\ \Omega_\ep ,\ \Grad v^{\pm}_{\ep,
\lambda} \cdot \vc{n}|_{\partial \Omega_\ep} = 0,
\]
satisfying
\[
\lim_{r \to \infty} r \left( \partial_r \pm {\rm
i}\sqrt{\lambda_\ep} \right) v^{\pm}_{\ep, \lambda_\ep} = 0,
\]
and
\bFormula{sa7}
\| v^{\pm}_{\ep, \lambda_\ep} \|_{L^2(B_{3R} \cap \Omega_\ep)} = 1.
\eF

Since
\[
\Delta v^{\pm}_{\ep, \lambda_\ep} = - \lambda_\ep v^{\pm}_{\ep,
\lambda_\ep} \ \mbox{in}\ B_{3R} \setminus B_R,
\]
the standard elliptic estimates yield
\[
v^{\pm}_{\ep, \lambda_\ep} \to v^{\pm}_\lambda \ \mbox{in}\ C^m(K) \
\mbox{for any}\ m \geq 0\ \mbox{and any compact}\ K \subset B_{3R}
\setminus \Ov{B}_R,
\]
in particular,
\[
v^{\pm}_{\ep, \lambda_\ep} \to v^{\pm}_{\lambda} \ \mbox{in}\
C^m(\partial B_{2R}) \ \mbox{for any}\ m \geq 0,
\]
which, together with formula (\ref{sa5}) yields
\[
v^{\pm}_{\ep, \lambda_\ep} \to v^{\pm}_\lambda \ \mbox{in}\ C^m(K) \
\mbox{for any}\ m \geq 0\ \mbox{and any compact}\ K \subset R^3
\setminus \Ov{B}_R,
\]
where $v^{\pm}_\lambda$ satisfy
\bFormula{sa8}
\Delta v^{\pm}_{\lambda} + \lambda v^{\pm}_{\lambda} = 0 \
\mbox{in}\ R^3 \setminus \Ov{B}_R, \ \lim_{r \to \infty} r \left(
\partial_r \pm {\rm i}\sqrt{\lambda} \right) v^{\pm}_{\lambda} = 0.
\eF

Finally, we claim that $v^{\pm}_\lambda$ satisfies
\bFormula{sa9}
\Delta v^{\pm}_\lambda + \lambda v^{\pm}_\lambda = 0 \ \mbox{in}\
\Omega \cap B_{3R}, \ \Grad v^{\pm}_\lambda \cdot \vc{n}|_{\partial
\Omega} = 0.
\eF
In order to see (\ref{sa9}) observe that
\bFormula{sa10}
\int_{\Omega_\ep} \left( \Grad v^{\pm}_{\ep, \lambda_\ep} \cdot
\Grad {\psi} - \lambda_\ep v^{\pm}_{\ep, \lambda_\ep} {\psi} +
\frac{\varphi_\ep {\psi} }{\| w^{\pm}_{\ep, \lambda_\ep}
\|_{L^2(B_{3R} \cap \Omega_\ep)} } \right) \ {\rm d}x = 0
\eF
for any $\psi \in \DC( \Ov{\Omega} \cap B_{3R} ) $. Since
$v^{\pm}_{\ep, \lambda_\ep}$ are bounded by (\ref{sa7}), we conclude
that
\bFormula{sa11}
\| v^{\pm}_{\ep, \lambda_\ep} \|_{W^{1,2}(B_{3R} \cap \Omega_\ep)}
\leq c,
\eF
uniformly for $\ep \to 0$.

As $\Omega_\ep$ possess the uniform extension property we may assume
that (\ref{sa11}) holds in $R^3$, and, by virtue of (\ref{prop1}),
we may pass to the limit in (\ref{sa10}) for each fixed $\psi$ to
obtain the desired conclusion (\ref{sa9}). Moreover, by the same
token, relation (\ref{sa7}) implies that
\bFormula{sa12}
\| v^{\pm}_\lambda \|_{L^2(B_{3R} \cap \Omega)} = 1.
\eF
However, relations (\ref{sa8}), (\ref{sa9}) imply $v^{\pm}_\lambda
\equiv 0$ in contrast with (\ref{sa12}). Thus we have shown
(\ref{sa6}).

Summing up the previous discussion we may infer that
\bFormula{sa13}
0 \leq  \left< \left( w^-_{\ep, \lambda} - w^+_{\ep, \lambda}
\right) ;\varphi \right>_{\Ome} \leq c(a,b, \varphi) \ \mbox{for
any}\ 0 < a < b,\ \varphi \in \DC(\Ome),
\eF
uniformly for $\ep \to 0$. Moreover, repeating the arguments of the
proof of (\ref{sa6}), we conclude that
\bFormula{sa14}
\lim_{\ep \to 0} \left< \left( w^-_{\ep, \lambda} - w^+_{\ep,
\lambda} \right) ;\varphi \right>_{\Ome} = \left< \left( w^-_{
\lambda} - w^+_{\lambda} \right) ;\varphi \right>_{\Omega} \
\mbox{for any}\ \lambda > 0,\ \varphi \in \DC(\Omega),
\eF
where
\[
\Delta w^{\pm}_{\lambda} + \lambda w^{\pm}_{\lambda} = \varphi \
\mbox{in}\ \Omega,\ \Grad w^{\pm}_{\lambda} \cdot \vc{n}|_{\partial
\Omega} = 0,
\]
\[
\lim_{r \to \infty} r \left( \partial_r \pm {\rm i} \sqrt{\lambda}
\right) w^{\pm}_{\lambda} = 0,\ r \equiv |x|.
\]

Seeing that
\bFormula{sa15}
\| \Grad \varphi \|_{L^2(\Ome;R^3)} = \| \sqrt{ - \Den } [\varphi]
\|_{L^2(\Ome; R^3)} ,
\eF
and
\[
\left< G(-\Den) [\varphi] ; \varphi \right>_{\Ome} = \int_0^\infty
G(\lambda) \ {\rm d} \mu_{\ep, \varphi }
\]
for any $\varphi \in \DC({\Omega}_\ep )$, we deduce from
(\ref{sa13}), (\ref{sa14}) that
\bFormula{sa16}
\| G(-\Den) [ \varphi ] \|_{W^{1,2}(\Ome)} \leq c(G, \varphi),
\eF
and
\bFormula{sa17}
1_{\Ome} G(-\Den)[\varphi] \to 1_{\Omega} G(-\Delta_N)[\varphi] \
\mbox{in}\ L^2(R^3),
\eF
for any $G \in \DC(0,\infty)$, $\varphi \in \DC(\Omega)$.

\subsection{Uniform decay for $|x| \to \infty$}

Our goal in this section is to establish the following decay
estimate:

\bLemma{a1} For $G \in \DC(0, \infty)$, $\varphi \in \DC(\Omega)$,
${\rm supp}[\varphi] \subset B_R$, we have
\[
\int_{|x| \geq R} |x|^{2s} |G(\sqrt{-\Den})[ \varphi ] |^2 \ {\rm
d}x \leq c(G,s,R) \| \varphi \|^2_{L^2(\Omega)} \ \mbox{for any}\ s
\geq 0,
\]
uniformly for $\ep \to 0$.
\eL

\bProof Extending $G$ as an \emph{even} function on $R$ we have
\[
G(\sqrt{-\Den})[\varphi] = \frac{1}{2} \int_{-\infty}^{\infty}
\widehat{G}(t) \left( \exp({\rm i} \sqrt{- \Den } t )  + \exp(-{\rm
i} \sqrt{- \Den } t )\right)[\varphi ] \ {\rm d}t,
\]
where $\widehat{G}$ stands for the Fourier transform of $G$.

Denoting
\[
\| v \|^2_{s,R} = \int_{|x| > R} |v|^2 |x|^{2s} \ {\rm d}x,
\]
we get
\[
\left\| G(\sqrt{-\Den})[\varphi] \right\|_{s,R}
\]
\[
\leq \frac{1}{2} \int_{-\infty}^{\infty}| \widehat{G}(t) | \left\|
\left( \exp({\rm i} \sqrt{- \Den } t )  + \exp(-{\rm i} \sqrt{- \Den
} t )\right)[\varphi ] \right\|_{s,R} {\rm d}t,
\]
where
\[
\left\| \left( \exp({\rm i} \sqrt{- \Den } t )  + \exp(-{\rm i}
\sqrt{- \Den } t )\right)[\varphi ] \right\|_{s,R}^2
\]
\[
= \intOe{ {\rm sgn}^+(|x| - R) |x|^{2s} \left| \left( \exp({\rm i}
\sqrt{- \Den } t )  + \exp(-{\rm i} \sqrt{- \Den } t
)\right)[\varphi ] \right|^2 }.
\]

On the other hand, the wave operator
\[
\left( \exp({\rm i} \sqrt{- \Den } t )  + \exp(-{\rm i} \sqrt{- \Den
} t )\right) = 2 \cos (\sqrt{ - \Den} t )
\]
admits a finite speed of propagation $1$, specifically,
\[
{\rm supp} \left[ \left( \exp({\rm i} \sqrt{- \Den } t )  +
\exp(-{\rm i} \sqrt{- \Den } t )\right) [\varphi] \right] \subset
B_{R + |t|};
\]
whence
\[
\intOe{ {\rm sgn}^+(|x| - R) |x|^{2s} \left| \left( \exp({\rm i}
\sqrt{- \Den } t )  + \exp(-{\rm i} \sqrt{- \Den } t
)\right)[\varphi ] \right|^2 }
\]
\[
\leq (|t| + R)^{2s} \intOe{  \left| \left( \exp({\rm i} \sqrt{- \Den
} t )  + \exp(-{\rm i} \sqrt{- \Den } t )\right)[\varphi ] \right|^2
}
\]
\[
= (|t| + R)^{2s} \| \varphi \|^2_{L^2(\Ome)}.
\]

As $G \in \DC(0,\infty)$, we have $(|t| + R)^{s} \widehat{G} \in
L^1(R)$ for any $s$, which completes the proof.

\qed

\subsection{Functional calculus - decay estimates}

Our ultimate goal in this section is the following result.

\bLemma{sa2}
We have
\[
\int_0^T \left| \left< \exp \left({\rm i} \sqrt{-\Den} \frac{t}{\ep}
\right)[\Psi], G(- \Den )[\varphi] \right>_{\Ome} \right|^2 \ \dt \leq \ep
c(\varphi,G) \| \Psi \|^2_{L^2(\Ome)}
\]
for any $\varphi \in \DC(\Omega)$, $\Psi \in L^2(\Ome)$, and any $G \in \DC(0,\infty)$.
\eL

\bProof We adapt the arguments of \cite{EF84}. By virtue \emph{spectral theorem}
(see Reed and Simon \cite[Chapter VIII]{ReSi1})
we have
\[
\left< \exp \left({\rm i} \sqrt{-\Den} \frac{t}{\ep} \right)[\Psi], G(-\Den) \varphi
\right>_{\Ome}
\]
\[
= \int_0^\infty \exp \left( {\rm i} \sqrt{\lambda}
\frac{t}{\ep} \right) G(\lambda) \tilde \Psi_\ep (\lambda) \ {\rm d} \mu_{\ep, \varphi}
(\lambda),
\]
where $\mu_{\ep,\varphi}$ is the spectral measure associated to the
function $\varphi$, and
\[
\tilde \Psi_\ep \in L^2(0, \infty; {\rm d}\mu_{\varphi}),\ \| \tilde
\Psi_\ep \|_{L^2_{\mu_\varphi}} \leq \| \Psi \|_{L^2(\Ome)}.
\]

Following Last \cite{Last} we deduce
\[
\int_0^T \left| \left< \exp \left({\rm i} \sqrt{-\Den} \frac{t}{\ep}
\right)[\Psi], G(-\Den) [\varphi] \right> \right|^2 \ \dt
\]
\[
\leq e T \sqrt{\pi} \int_0^\infty \int_0^\infty \tilde \Psi_\ep (x) \
\Ov{ \tilde \Psi_\ep (y) } \exp \left( - \frac{ T^2 | \sqrt{x} -
\sqrt{y} |^2 }{4 \ep^2} \right) G(x) G(y) \ {\rm d}\mu_{\ep,\varphi} (x) \
{\rm d}\mu_{\ep, \varphi} (y);
\]
therefore, by Cauchy-Schwartz inequality,
\[
\int_0^T \left| \left< \exp \left({\rm i} \sqrt{-\Den} \frac{t}{\ep}
\right)[\Psi], G(-\Den)[\varphi] \right>_{\Ome} \right|^2 \ \dt
\]
\[
\leq e T \sqrt{\pi} \int_0^\infty |\Psi_\ep (x) |^2 \left( \int_0^\infty
\exp \left( - \frac{ |\sqrt{x} - \sqrt{y}|^2 }{\ep^2} \frac{T^2}{4}
\right) {\rm d}\mu_{\ep,\varphi} (y) \right)  G^2(x) \ {\rm d}\mu_{\ep,\varphi}
(x).
\]

Furthermore,
\bFormula{sa18}
\int_0^\infty \exp \left( - \frac{ |\sqrt{x} - \sqrt{y}|^2 }{\ep^2}
\frac{T^2}{4} \right) {\rm d}\mu_{\ep,\varphi} (y)
\eF
\[
= \sum_{n=0}^\infty
\int_{ \ep n \leq |\sqrt{y} - \sqrt{x} | < \ep(n+1) }  \exp \left( -
\frac{ |\sqrt{x} - \sqrt{y}|^2 }{\ep^2} \frac{T^2}{4} \right) {\rm
d}\mu_{\ep,\varphi} (y)
\]
\[
\leq \sup_{n \geq 0} \int_{ \ep n \leq |\sqrt{y} - \sqrt{x} | <
\ep(n+1) } 1 {\rm d}\mu_{\ep,\varphi} (y) \sum_{n=0}^\infty \exp \left( -
\frac{n^2 T^2}{4} \right).
\]
Since only the points $x$ belonging to ${\rm supp}[G]$ are relevant
in (\ref{sa18}), the length of the intervals
\[
I_n(x) = \{ y \in [0, \infty) \ | \ \ep n \leq |\sqrt{y} - \sqrt{x}
| < \ep(n+1) \}
\]
never exceeds $\ep$,
\[
|I_n(x)| \leq c_G \ep.
\]
Thus we conclude combining (\ref{sa18}), with the uniform bounds on
the spectral measures (see (\ref{sa2})) established in (\ref{sa13}).

\qed

\section{Dispersive estimates, local decay of acoustic waves}
\label{co}

Returning to the acoustic equation (\ref{ae19}), (\ref{ae20}) we can use the dispersive estimate established in Lemma \ref{Lsa2} to show that
\bFormula{dec1}
\left\{ t \mapsto \intOe{ \Phi_\ep G(-\Den) [\varphi ] } \right\} \to 0
\ \mbox{in}\ L^2(0,T)
\eF
for any fixed $G \in \DC(0, \infty)$, $\varphi \in \DC(\Omega)$. Indeed, by means of Duhamel's formula, we have
\[
\Phi_\ep(t,\cdot) = \frac{1}{2} \exp \left({\rm i} \sqrt{-\Den} \frac{t}{\ep}
\right) \left[ \Phi_{0,\ep} + {\rm i} \frac{1}{\sqrt{- \Den}} [r_{0,\ep}]
\right]
\]
\[
+ \frac{1}{2} \exp \left(-{\rm i} \sqrt{-\Den} \frac{t}{\ep} \right)
\left[ \Phi_{0,\ep} - {\rm i} \frac{1}{\sqrt{- \Den}} [r_{0,\ep}]
\right]
\]
\[
+ \ep^{-2\beta}
\frac{1}{2} \int_0^t \left( \exp \left({\rm i} \sqrt{-\Den} \frac{t - s}{\ep}
\right) + \exp \left(-{\rm i} \sqrt{-\Den} \frac{t - s}{\ep}
\right)\right)[H_\ep(s)] \ {\rm d}s,
\]
with
\[
H_\ep = G^1_\ep + (-\Den)^{-1/2}[ G^2_\ep] + (-\Den)^{1/2} [G^3_\ep]
+ (-\Den)[G^4_\ep] + (-\Delta_{\ep,N})^{-1} [G^5_\ep],
\]
(see (\ref{ae20})), where we have assumed, for the sake of
simplicity, that $p'(\Ov{\vr}) = 1$. Consequently, the desired
conclusion (\ref{dec1}) follows from  Lemma \ref{Lsa2} as soon as
$\beta < 1/4$, see \cite{EF84} for details.

\subsection{Compactness in time of the momenta}

In view of estimate (\ref{ub6}), the desired strong convergence
claimed in (\ref{mm3}) follows as soon as we show that
\bFormula{ct1}
\left\{ t \mapsto \intO{ \vc{V}_\ep(\cdot, t) \cdot \vc{w} } \right\}
\to \left\{ t \mapsto \Ov{\vr} \intO{ \vc{U}(\cdot,t) \cdot \vc{w} } \right\} \ \mbox{in} \ L^2(0,T) \ \mbox{for}\
\vc{w} \in \DC(\Omega;R^3).
\eF
Indeed relation (\ref{ct1}) with (\ref{convv}) imply that
\[
\int_0^T \int_K \vre |\vue|^2 \ \dxdt \to \Ov{\vr} \int_0^T \int_K
|\vc{U} |^2 \ \dx \ \dt \ \mbox{for any compact}\ K \subset \Omega,
\]
yielding (\ref{mm3}).

In order to show (\ref{ct1}), we use Helmholtz decomposition to obtain
\[
\intOe{ \vc{V}_\ep \cdot \vc{w} } = \intOe{ \vc{H}_\ep [\vre \vue] \cdot \vc{w} } -
\intOe{ \Phi_\ep \Div \vc{w} }
\]
\[
=\intOe{ \vre \vue \cdot \vc{H}_\ep [\vc{w}] } -
\intOe{ \Phi_\ep \Div \vc{w} }
\]
\[
= \intOe{ \vre \vue \cdot \vc{H} [\vc{w}] }
+ \intOe{ \vre \vue \cdot (\vc{H}_\ep[\vc{w}] - \vc{H} [\vc{w}]) } - \intOe{ \Phi_\ep \Div \vc{w} },
\]
where, in accordance with (\ref{ae2}) and the standard Aubin-Lions
argument,
\bFormula{ct2}
\left\{ t \mapsto \intOe{ \vre \vue \cdot \vc{H} [\vc{w}] } \right\} \to
\left\{ t \mapsto \Ov{\vr}\intO{ \vc{U} \cdot \vc{w} } \right\}\ \mbox{in} \ L^2(0,T).
\eF
Here, we have extended $\vc{H}[\vc{w}]$ and $\vc{H}_\ep [\vc{w}]$ by zero outside
$\Omega$ and $\Ome$, respectively.

Furthermore, we write
\[
\intOe{ \vre \vue \cdot (\vc{H}_\ep[\vc{w}] - \vc{H} [\vc{w}]) }
\]
\[
=
\intOe{ (\vre - \Ov{\vr})  \vue \cdot (\vc{H}_\ep[\vc{w}] - \vc{H} [\vc{w}]) }
+ \Ov{\vr} \intOe{ \vue \cdot (\vc{H}_\ep[\vc{w}] - \vc{H} [\vc{w}]) } .
\]
In view of estimates (\ref{ae5}), (\ref{ae5a}), and boundedness of
$\vc{H}_\ep$, $\vc{H}$  in $L^p$ (see (\ref{ee4}), (\ref{ee5})), we
get
\bFormula{ct3}
\left\{ t \mapsto  \intOe{ (\vre - \Ov{\vr})  \vue \cdot (\vc{H}_\ep[\vc{w}] - \vc{H} [\vc{w}]) } \right\} \to 0 \ \mbox{in}\ L^2(0,T).
\eF
Moreover, as
\[
\vc{H}_\ep [\vc{w} ] \to \vc{H} [\vc{w}] \ \mbox{weakly in} \ L^2(\Omega;R^3),
\]
we get, by virtue of (\ref{mm2}),
\bFormula{ct4}
\left\{ t \mapsto \intOe{ \vue \cdot (\vc{H}_\ep[\vc{w}] - \vc{H} [\vc{w}]) } \right\} \to 0 \ \mbox{in}\ L^2(0,T).
\eF

Finally,
\[
\intOe{ \Phi_\ep \Div \vc{w} } = \intOe{ \Phi_\ep G(-\Den) [ \Div \vc{w} ] } +
\intOe{ \Phi_\ep (1 - G(-\Den))[ \Div \vc{w} ] },
\]
where, as stated in (\ref{dec1}),
\bFormula{ct5}
 \left\{ t \mapsto \intOe{ \Phi_\ep G(-\Den) [ \Div \vc{w} ] } \right\} \to 0
 \ \mbox{in}\ L^2(0,T).
\eF

Writing $\Phi_\ep = \Phi^1_\ep + \Phi^2_\ep$ as in (\ref{ae9}), (\ref{ae10}), we have
\bFormula{ct6}
\left\{ t \mapsto \intOe{ \Phi^2_\ep (1 - G(-\Den))[ \Div \vc{w} ] } \right\} \to 0
\ \mbox{in} \ L^2(0,T) ,
\eF
while, in agreement with (\ref{sa17}) and Lemma \ref{La1},
\bFormula{ct7}
\left\{ t \mapsto \intOe{ \Phi^1_\ep (1 - G(-\Den))[ \Div \vc{w} ] }
\right\}
\eF
\[ \to
\left\{ t \mapsto \intO{ \Phi^1 (1 - G(-\Delta_N))[ \Div \vc{w} ] }
\right\} \ \mbox{in} \ L^2(0,T) ,
\]
where the resulting expression is small as soon as $G \approx 1_{[0,
\infty)}$. Indeed
\[
\intO{\Phi^1 (1 - G(-\Delta_N))[ \Div \vc{w} ] } = \intO{
(-\Delta_N)^{1/2} \Phi^1 \frac{1}{(-\Delta_N)^{1/2} }(1 -
G(-\Delta_N))[ \Div \vc{w} ] },
\]
where
\[
\| (-\Delta_N)^{1/2} [\Phi^1] \|_{L^2(\Omega)} = \| \Grad \Phi^1
\|_{L^2(\Omega)},
\]
while
\[
\left\| \frac{1}{(-\Delta_N)^{1/2}} [\Div \vc{w} ]
\right\|_{L^2(\Omega)} = \left\| {(-\Delta_N)^{1/2}}
(-\Delta_N)^{-1} [\Div \vc{w} ] \right\|_{L^2(\Omega)}
\]
\[
= \left\| \Grad (-\Delta_N)^{-1} [\Div \vc{w} ]
\right\|_{L^2(\Omega)} \leq \| \vc{w} \|_{L^2(\Omega;R^3)}.
\]

Relations (\ref{ct2} - \ref{ct7}) imply (\ref{ct1}), in particular,
we have shown (\ref{mm3}). The proof of Theorem \ref{Tm1} is now
complete.

\subsection{Boundary behavior of the limit velocity field $\vc{U}$}

\label{coa}

In the previous analysis, we left open the problem of the boundary
conditions satisfied by the limit velocity field $\vc{U}$. To this
end, we revoke the results of \cite{BFNW}. Suppose that, after
suitable translation and rotation of the coordinate system, a part
$\Gamma$ of the boundary of $\Omega$ can be described by a graph of
a function $b \in W^{1,\infty} (U)$, $U \subset R^2$,
\[
\Gamma = \{ (x_1,x_2,x_3) \ | \ (x_1,x_2) \in U, \ x_3 = b(x_1,x_2) \},
\]
while $\Gamma_\ep = \partial \Ome \cap U \times R$ are represented as
\[
\Gamma_\ep = \{ (x_1,x_2,x_3) \ | \ (x_1,x_2) \in U, \ x_3 = b_\ep(x_1,x_2) \},
\]
where $\{ b_\ep \}_{\ep>0} $ is a bounded sequence in $W^{1,\infty}(U)$, $b_\ep \to b$ in
$C(\Ov{U})$. Similarly to \cite{BFNW}, we assume that the boundaries $\Gamma_\ep$
are oscillating for $\ep \to 0$. More specifically, introducing a \emph{Young measure}
${\cal R}[y]$, $y \in U$, associated to the gradients $\{ \nabla_y b_\ep\}_{\ep > 0}$, we suppose that
\bFormula{young}
{\rm supp}[ {\cal R}[y]] \ \mbox{contains two independent vectors in} \ R^2
\ \mbox{for a.a.} \ y \in U.
\eF
As shown in \cite{BFNW}, condition (\ref{young}) implies
\[
\vc{U}|_{\Gamma} = 0,
\]
see also B\v rezina \cite{Brez} for refined results in this
direction.

\end{document}